\documentclass[10pt]{amsart}

\usepackage{graphicx}
\usepackage{latexsym,amsmath}
\usepackage{amsmath,amsthm}
\usepackage{amsfonts}
\usepackage[psamsfonts]{amssymb}
\usepackage{enumerate}

\usepackage{url}

\usepackage{pstricks}

\theoremstyle{plain}
\newtheorem{theorem}{Theorem}[section]
\newtheorem{corollary}[theorem]{Corollary}

\newtheorem{lemma}[theorem]{Lemma}
\newtheorem{proposition}[theorem]{Proposition}

\theoremstyle{definition}

\theoremstyle{definition}

\theoremstyle{remark}

\theoremstyle{remark}

\newtheorem*{remark*}{Remark}
\numberwithin{equation}{section}

\newcommand{\R}{\mathbb{R}}
\newcommand{\Q}{\mathbb{Q}}
\renewcommand{\P}{\mathsf{P}}
\newcommand{\si}{\sigma}
\newcommand{\ka}{\kappa}
\newcommand{\la}{\lambda}
\newcommand{\vp}{\varepsilon}
\newcommand{\ts}{\tilde S}

\renewcommand{\le}{\leqslant}

\newcommand{\circl}{\,%
\includegraphics[bb=
0 0 595 842,width=7pt,height=8pt]{circl.eps}}

\begin{document}

%
%

\title[Non-degeneracy of correlation coefficients]
{
On the Non-degeneracy of Kendall's and Spearman's Correlation Coefficients}
\date{\today}
\author{Iosif Pinelis}
\address{
Department of Mathematical Sciences\\
Michigan Technological University\\
Houghton, Michigan 49931 }
\email{ipinelis@mtu.edu}
\keywords{Kendall's correlation coefficient, Spearman's correlation coefficient, asymptotic variance, non-degeneracy, support of distribution}


\begin{abstract}
Hoeffding proved that Kendall's and Spearman's nonparametric measures of correlation between two continuous random variables $X$ and $Y$ are each asymptotically normal with an asymptotic variance of the form $\sigma^2/n$ -- provided the non-degeneracy condition $\sigma^2>0$ holds, where $\sigma^2$ is a certain (always nonnegative) expression which is determined by the joint distribution (say $\mu$) of $X$ and $Y$. Sufficient conditions for $\sigma^2>0$ in terms of the support set (say $S$) of $\mu$ are given, the same for both correlation statistics. One of them is that there exist a rectangle with all its vertices in $S$, sides parallel to the $X$ and $Y$ axes, and an interior point also in $S$. Another sufficient condition is that the Lebesgue measure of $S$ be nonzero. 
\end{abstract}

\maketitle

\section{Introduction}\label{intro}
Let $(X,Y)$ be a random point in $\R^2$ with a (joint) cumulative distribution function (c.d.f.) $F$ and continuous marginal c.d.f.'s $F_X$ and $F_Y$, so that $F(x,y)=\P(X\le x,Y\le y)$, $F_X(x)=\P(X\le x)$, and $F_Y(y)=\P(Y\le y)$ for all real $x$ and $y$. Then $F$ is continuous as well, since
$|F(x_2,y_2)-F(x_1,y_1)|\le |F_X(x_2)-F_X(x_1)|+ |F_Y(y_2)-F_Y(y_1)|$
for all $x_1,y_1,x_2,y_2$ in $\R$.
Vice versa, if $F(x,y)$ is continuous in $x$ for each real $y$, then $0=F(x,y)-F(x-,y)=\P(X=x,Y\le y)\underset{y\to\infty}\longrightarrow
\P(X=x)=F_X(x)-F_X(x-)$, so that $F_X$ is continuous. Similarly, if $F(x,y)$ is continuous in $y$ for each real $x$, then $F_Y$ is continuous. So, the marginal c.d.f.'s $F_X$ and $F_Y$ are continuous iff the joint c.d.f.\ $F$ is so; in such a case, one may simply say that the distribution of $(X,Y)$ is continuous.

Let $\mu=\mu_{X,Y}$ denote the measure that is the probability distribution of $(X,Y)$.
Let $S=S_{X,Y}$ stand for the set in $\R^2$ that is the support of $\mu$ \big(defined as the intersection of all closed sets of $\mu$-measure 1; then $S$ is the smallest of all such sets, and also $S$ coincides with the set of all $x\in\R^2$ such that  $\mu(B_\vp(x))>0$ for all $\vp>0$, where $B_\vp(x)$ is the (say open) disk of radius $\vp$ centered at $x$\big). 

The most common nonparametric statistics that measure association between $X$ and $Y$ are Spearman's rank correlation \cite{sp} and Kendall's difference sign correlation \cite{ken}, usually denoted by $\rho$ and $\tau$, respectively. These statistics are based on a sample of
independent random points $(X_1,Y_1),\dots,(X_n,Y_n)$ each having the same distribution as the random point $(X,Y)$. 
%
In his landmark paper \cite{ho}, Hoeffding proved that $\rho$ and $\tau$ are each asymptotically normal as $n\to\infty$ with an asymptotic variance of the form $\si_\ka^2/n$ --- provided the non-degeneracy condition $\si_\ka^2>0$ holds, where $\ka$ is either $\rho$ or $\tau$, and 
$\si_\ka^2$ is a certain (always nonnegative) expression which is determined by the c.d.f.\ $F$.
It is therefore important to have convenient criteria for the non-degeneracy condition $\si_\ka^2>0$. 

\section{Kendall's $\tau$} 
As follows from
Hoeffding \cite[(9.13)]{ho}, $\si^2_\tau=0$ iff the function $d_\tau$ defined by the formula
$$d_\tau(x,y):=F(x,y)-\big(F_X(x)+F_Y(y)\big)/2$$
is constant on a set (say $A$) of $\mu$-measure 1; then $d_\tau$ must be constant on the closure of $A$ (since $d_\tau$ is a continuous function) and hence on $S$.
That is, one has
\begin{proposition}\label{prop:}
$\si_\tau^2=0$ iff $d_\tau$ is constant on $S$. 
\end{proposition}

Even in such simple cases as the ones considered in examples given at the end of this section, it may be not quite immediately obvious based on Proposition~\ref{prop:}
whether the distribution of $\tau$ is asymptotically non-degenerate, in the sense that $\si_\tau^2>0$.
However, we shall give simple conditions sufficient for such non-degeneracy.

In what follows, the term \emph{rectangle} means a nonempty set of the form $(x_1,x_2]\times(y_1,y_2]$; the points $(x_1,y_1)$, $(x_1,y_2)$, $(x_2,y_1)$, $(x_2,y_2)$ in $\R^2$ are then naturally called the vertices of the rectangle.

\begin{lemma}\label{th:}
Suppose that $S$ contains all the four vertices and also an interior point of a rectangle; then $\si_\tau^2>0$. 
\end{lemma}

\begin{proof}
Assume, to the contrary, that $\si_\tau^2=0$. Then, by Proposition~\ref{prop:}, $d_\tau$ is constant on $S$.
On the other hand, one has 5 points $(x_*,y_*)$, $(x_1,y_1)$, $(x_1,y_2)$, $(x_2,y_1)$, $(x_2,y_2)$ in $S$ such that $x_1<x_*<x_2$ and $y_1<y_*<y_2$. So, in view of the definition of $S$, 
\begin{equation}\label{eq:main}
\begin{aligned}
0<\mu\big((x_1,x_2)\times(y_1,y_2)\big)
&\le\mu\big((x_1,x_2]\times(y_1,y_2]\big) \\
&=d_\tau(x_2,y_2)-d_\tau(x_1,y_2)-d_\tau(x_2,y_1)+d_\tau(x_1,y_1)=0,
\end{aligned}	
\end{equation}
which is a contradiction.
\end{proof}

One simple sufficient condition is an immediate corollary of Lemma~\ref{th:}:
\begin{corollary}\label{cor:0}
If the interior of $S$ is non-empty, then $\si_\tau^2>0$.
\end{corollary}


Working a bit harder, one can get a stronger result.

Let $\la_k$ denote the Lebesgue measure for $\R^k$.
\begin{corollary}\label{cor:1}
Suppose that $\la_2(S)>0$; then $\si_\tau^2>0$.
\end{corollary}
\begin{proof}
Since $\R^2$ can be partitioned into countable many disjoint rectangles, one has $\la_2(\ts)>0$ for some rectangle $R$, where $\ts:=R\cap S$. Further, $\ts$ can be approximated by the union of disjoint rectangles $R_1,R_2,\dots$ contained in $R$ so that $\ts\subseteq\bigcup_n\,R_n$ and
$\la_2(\ts)>\frac23\,\la_2\big(\bigcup_n\,R_n\big)$, that is, 
$\sum_n\big(\la_2(\ts\cap R_n)-\frac23\la_2(R_n)\big)>0$, whence
$\la_2(\ts\cap R_m)>\frac23\la_2(R_m)$ for some natural $m$. Shifting and re-scaling if necessary, without loss of generality let us assume that $R_m=(0,1]^2:=(0,1]\times(0,1]$. Then, introducing the set $M:=\ts\cap R_m$, one has
$$M\subseteq S\cap(0,1]^2\quad\text{and}\quad\la_2(M)>\tfrac23.$$
Now, by Fubini's theorem,
$\int_0^1\la_1(M^x)\,dx=\la_2(M)>\frac23$, where $M^x:=
\{y\in\R\colon(x,y)\in M\}$.

Hence, the set $A:=\{x\in(0,1)\colon\la_1(M^x)>\frac23\}$ is infinite \big(otherwise, one would have $\la_1(M^x)\le\frac23$ for almost all $x$ in $(0,1)$ and thus $\int_0^1\la_1(M^x)\,dx\le\frac23$\big).
Therefore, there are $x_1$, $x_*$, $x_2$ in $A$ such that $x_1<x_*<x_2$. Then, for $M^*:=M^{x_1}\cap M^{x_*}\cap M^{x_2}$, one has $\la_1(M^*)>1-3(1-\frac23)=0$, so that the set $M^*$ is infinite and thus contains some $y_1$, $y_*$, $y_2$ such that $y_1<y_*<y_2$. It remains to refer to Lemma~\ref{th:}.
\end{proof}

\begin{corollary}\label{cor:2}
Suppose that the measure $\mu$ has a nonzero absolutely continuous component (with respect to the Lebesgue measure $\la_2$); then $\si_\tau^2>0$.
\end{corollary}

\begin{proof}
Let $f$ be the density of the nonzero absolutely continuous component of $\mu$. Then $0<\int_{\R^2}f\,d\la_2=\int_S f\,d\la_2$ \big(since $0\le\int_{\R^2\setminus S}f\,d\la_2\le\mu(\R^2\setminus S)=0$\big). Therefore, $\la_2(S)>0$. It remains to refer to Corollary~\ref{cor:1}.
\end{proof}

Observe that Corollary~\ref{cor:1} \emph{strictly} contains Corollary~\ref{cor:2}. Indeed, there is a random point $(X,Y)$ with a continuous c.d.f.\ $F$ such that $\la_2(S)>0$ while $\mu$ is singular with respect to $\la_2$. For example, let $X$ be any random variable with an everywhere strictly positive density (with respect to $\la_1$). Next, let $Y:=X+Q$, where $Q$ is any random variable with values in the set $\Q$ of all rational real numbers such that $\P(Y=r)>0$ for all  $r\in\Q$. Then the random variable $Y$ is absolutely continuous and $\P((X,Y)\in S_0)=1$, where $S_0:=\bigcup_{r\in\Q}\{(x,x+r)\colon x\in\R\}$. At that, $\la_2(S_0)=0$, so that the measure $\mu$ is singular with respect to $\la_2$, whereas $S=\R^2$ and hence $\la_2(S)>0$.

This example also shows that Corollary~\ref{cor:2} does not even contain Corollary~\ref{cor:0}. On the other hand, it is easy to see that, vice versa, Corollary~\ref{cor:0} does not contain Corollary~\ref{cor:2}; indeed, let $(X,Y)$ be uniformly distributed on a set of the form $C\times C$, where $C$ is any ``fat'' Cantor subset of $\R$ (e.g.\ the so-called Smith-Volterra-Cantor set, \url{http://en.wikipedia.org/wiki/Smith-Volterra-Cantor_set}), which is a non-empty compact nowhere-dense set such that $\la_1(C\cap B_\vp(x))>0$ for all $x\in C$ and all $\vp>0$;
then the measure $\mu$ is absolutely continuous with respect to $\la_2$, whereas $S=C\times C$, so that the interior of $S$ is empty. The latter example also shows that Corollary~\ref{cor:1} \emph{strictly} contains Corollary~\ref{cor:0}.

Note that the 5-point condition in Lemma~\ref{th:} is not necessary for $\si_\tau^2>0$. For example, suppose that the support of the measure $\mu$ is the union of straight line segments $S_1:=\{(x,x)\colon0\le x\le1\}$ and $S_2:=\{(x,1-x)\colon\frac12\le x\le1\}$. Then the 5-point condition in Lemma~\ref{th:} is not satisfied, and yet, the function $d_\tau$ is not constant on either segment $S_1$ or $S_2$, so that, by Proposition~\ref{prop:}, $\si_\tau^2>0$.

However, the role of the interior point in the \emph{5}-point condition in Lemma~\ref{th:} is in a certain sense indispensable.

\vspace*{8pt}





\begin{center}
{
\psset{unit=4}%
\begin{pspicture}(0,0)(1,1)
\psgrid[subgriddiv=0,griddots=40,gridlabels=0pt]
\psset{unit=10mm}
\psset{origin={0,0}}
\psset{linewidth=.2mm}
\psset{linewidth=.2mm}
\psline[linewidth=.3mm](2,0)(4,2)
\psline[linewidth=.3mm](4,2)(2,4)
\psline[linewidth=.3mm](2,4)(0,2)
\psline[linewidth=.3mm](0,2)(2,0)
\psline[linestyle=dashed,dash=.1 .1](.7,1.3)(3.3,1.3)
\psline[linestyle=dashed,dash=.1 .1](3.3,1.3)(3.3,2.7)
\psline[linestyle=dashed,dash=.1 .1](3.3,2.7)(.7,2.7)
\psline[linestyle=dashed,dash=.1 .1](.7,2.7)(.7,1.3)
\uput[270](2,-.1){(a)}
\end{pspicture}
}
\hspace{1.5cm}
\begin{pspicture}
(0,0)(4,4)
\psgrid[subgriddiv=0,griddots=10,gridlabels=0pt]
\psset{unit=10mm}
\psset{linewidth=.2mm}
\pscircle[](.5,.5){.5}
\psarc[](2,3){1}{180}{270}
\psarc[](4,1){1}{90}{180}
\psarc[](2,4){.5}{270}{360}
\psarc[](3,4){.5}{180}{270}
\psarc[](2,3){.5}{0}{90}
\psarc[](3,3){.5}{90}{180}
\psset{linewidth=.2mm}
\uput[270](2,-.1){(b)}
\end{pspicture}

\vspace*{21pt}

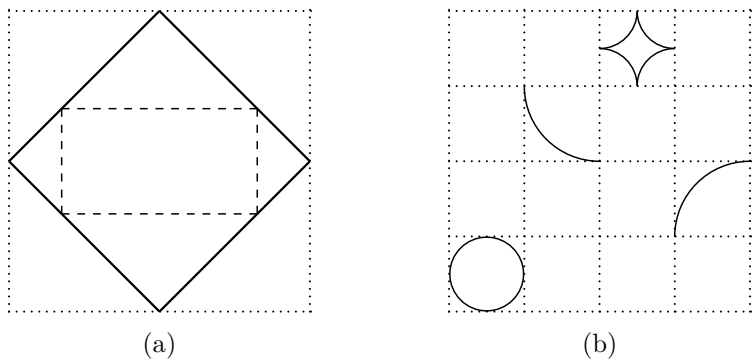
\textsc{Figure 1.} Examples with $d_\tau$ constant on $S$, 
that is, with $\si_\tau^2=0$.
\end{center}

\vspace*{5pt}

For example, let the random point $(X,Y)$ be uniformly distributed on the union $S$ of the four sides of the thick-line square shown in Fig.~1 (a). (More generally, one may assume that the distribution is continuous, symmetric about
the 
center of
the thick-line square, and has $S$ as its support.)
Obviously, here there are infinitely many rectangles with all the four vertices in $S$; one of them is shown here (dashed); however, there is no 5th, interior point in any such rectangle which would also belong to the support $S$ of the distribution of $(X,Y)$.
Accordingly, here $\si_\tau^2=0$ by Proposition~\ref{prop:}, since $d_\tau(X,Y)=-\frac14$ with probability 1.

\renewcommand{\u}{\,%
\includegraphics[bb=0 0 9 10,width=7pt,height=8pt]{u.eps}}

\renewcommand{\circl}{\rule{0pt}{13pt}%
\begin{picture}(10,10)
\linethickness{.5mm}
\put(5,5) {\circle{10}}
\end{picture}
}

\renewcommand{\circl}{\begin{pspicture}(0,.05)(.25,.4)
\psset{unit=10pt}
\psset{linewidth=.1mm}
\pscircle[](.5,.5){.5}
\end{pspicture}
}

\renewcommand{\star}{\begin{pspicture}(0,.05)(.25,.4)
\psset{unit=10pt}
\psset{linewidth=.1mm}
\psarc[](0,1){.5}{270}{360}
\psarc[](1,1){.5}{180}{270}
\psarc[](0,0){.5}{0}{90}
\psarc[](1,0){.5}{90}{180}
\end{pspicture}
}

\newcommand{\down}{\begin{pspicture}(0,.05)(.25,.4)
\psset{unit=10pt}
\psset{linewidth=.1mm}
\psarc[](1,1){1}{180}{270}
\end{pspicture}
}

\newcommand{\up}{\begin{pspicture}(0,.05)(.25,.4)
\psset{unit=10pt}
\psset{linewidth=.1mm}
\psarc[](1,0){1}{90}{180}
\end{pspicture}
}

Another example is illustrated by Fig.~1 (b), where it is assumed that the support of the distribution of the random point $(X,Y)$ is the union of all the four thick-line components: \circl, \down, \star, and \up, whose probability masses are $\frac6{11}$, $\frac1{11}$, $\frac2{11}$, and $\frac2{11}$, respectively, and at that the (conditional) distribution on either one of the pieces \circl\ and \star\ is continuous and symmetric about 
the center of the piece. The (conditional) distribution on either one of the pieces \down\ and \up\ does not matter at all, as long as it is continuous. Here too there are infinitely many rectangles with all the four vertices in $S$, but the interior of none of them has a nonzero $\mu$-mass.
Accordingly, $\si_\tau^2=0$, since $d_\tau(X,Y)=-\frac3{22}$ with probability 1. However, any weights of the pieces \circl, \down, \star, and \up other than $\frac6{11}$, $\frac1{11}$, $\frac2{11}$, and $\frac2{11}$ would result in $\si_\tau^2>0$.

These examples suggest that it would be difficult, if at all possible, to find a necessary and sufficient condition for the non-degeneracy from which one could easily deduce such results as Corollary~\ref{cor:1}. 

\section{Spearman's $\rho$}
By Hoeffding \cite[(9.27)]{ho}, $\si^2_\rho=0$ iff a certain function $d_\rho$ of the form 
$$d_\rho(x,y):=F_X(x)F_Y(y)-f(x)-g(y)$$
for some continuous functions $f$ and $g$ 
is constant on a set of $\mu$-measure 1 or, equivalently, on the support $S$. 

So, $\si_\rho^2=0$ implies that --- again for any 
5 points $(x_*,y_*)$, $(x_1,y_1)$, $(x_1,y_2)$, $(x_2,y_1)$, $(x_2,y_2)$ in $S$ such that $x_1<x_*<x_2$ and $y_1<y_*<y_2$ --- all the relations \eqref{eq:main} hold  
with $\mu_X\otimes\mu_Y$ and $d_\rho$ instead of $\mu$ and $d_\tau$, respectively; here, as usual, $\mu_X$ and $\mu_Y$ stand for the distributions of $X$ and $Y$, so that $\mu_X\otimes\mu_Y$ is the corresponding product measure. 

Thus, obtains the following result.

\begin{corollary}\label{cor:3}
Lemma~\ref{th:} holds with $\si_\rho$ in place of $\si_\tau$, and so do Corollaries~\ref{cor:0}--\ref{cor:2}. 
\end{corollary}

%
%
%
%
%
%


%
%
%
%
\end{document}